# Applications of Finite Markov Chain Models to Management[*]

**Michael Gr. Voskoglou**
Professor Emeritus of Mathematical Sciences
Graduate Technological Educational Institute of Western Greece
School of Technological Applications
Meg. Alexandrou 1 – 263.34 Patras- Greece
E - mail: voskoglou@teiwest.gr , mvosk@hol.gr

## Abstract

Markov chains offer ideal conditions for the study and mathematical modelling of a certain kind of situations depending on random variables. The basic concepts of the corresponding theory were introduced by Markov in 1907 on coding literary texts. Since then, the Markov chain theory was developed by a number of leading mathematicians, such as Kolmogorov, Feller etc. However, only from the 60's the importance of this theory to the Natural, Social and most of the other Applied Sciences has been recognized. In this review paper we present applications of finite Markov chains to Management problems, which can be solved, as most of the problems concerning applications of Markov chains in general do, by distinguishing between two types of such chains, the ergodic and the absorbing ones.

## 1. Introduction

One of the most interesting approaches of the nature and position of mathematics into the whole body of the human knowledge determines it as an activity whose target is to interpret and (or) describe the various phenomena of the real world, as well as the real situations of our everyday life. This approach, combined with the recent waves of social, economical, and technological changes and evolutions of our world, transubstantiates mathematics to an essential factor for the formation of the future of our society.

The above approach is materialized through the use of *mathematical models*, which are simplified representations of the corresponding real phenomena or situations, achieved through the use of mathematical terms and symbols, i.e. functions, equations, inequalities, etc. A mathematical model excels compared with the other types of models in use, e.g. simulation, analogical, iconic, etc [9], because it gives accurate and general solutions to the corresponding problems. Davis and Hersh [2] emphasize this fact by stating that: 'The usefulness of a mathematical model is exactly its success to foresee, or (and) to imitate accurately the behaviour of the real world".

---

[*] The contents of this paper were presented in an author's Keynote Speech at the *International Conference on Knowledge Engineering and Big Data Analytics (KE &BDA)* in Future University, Cairo, Egypt, December 15-16, 2015

A very important theory that offers ideal conditions for the study and mathematical modelling of a certain kind of phenomena depending on random variables is the theory of Markov Chains. Roughly speaking, a *Markov Chain* is a stochastic process that moves in a sequence of steps (phases) through a set of states and has a "one-step memory", i.e. the probability of entering a certain state in a certain step, although in practice may not be completely independent of previous steps, depends at most on the state occupied in the previous step. This is known as the *Markov property*. When the set of its states is a finite set, then we speak about a *finite Markov Chain*. For special facts on finite Markov Chains we refer to the classical on the subject book of Kemeny & Snell [4].

The basic concepts of Markov Chains were introduced by A. Markov in 1907 on coding literary texts. Since then the Markov Chains theory was developed by a number of leading mathematicians, such as A. Kolmogorov, W. Feller etc. However, only from the 1960's the importance of this theory to the Natural, Social and most of the other Applied Sciences has been recognized (e.g. see [1, 3, 8], etc). Most of the problems concerning applications of Markov Chains can be solved by distinguishing between two types of such Chains, the *Ergodic* and the *Absorbing* ones.

This review paper reports characteristic applications of finite Markov Chains to Management problems. The rest of the paper is formulated as follows: In Section 2 we introduce the basic form of the Markov Chain model, which is used to make short run forecasts for the evolution of various phenomena. In Section 3 we present the Ergodic Markov Chains, the equilibrium situation of which is used to obtain long run forecasts. Section 4 is devoted to applications of Absorbing Markov Chains to Management, while in Section 5 we give a brief account of our relevant researches on applications of Markov Chains to Education and Artificial Intelligence. Finally, in Section 6 we state our conclusions and we discuss the perspectives of future research on the subject.

## 2. The Basic Form of the Markov Chain Model

Let us consider a finite Markov Chain with n states, where n is a non negative integer, $n \geq 2$. Denote by $p_{ij}$ the *transition probability* from state $s_i$ to state $s_j$, i, j = 1, 2,…, n ; then the matrix **A= [$p_{ij}$]** is called the *transition matrix* of the Chain. Since the transition from a state to some other state (including itself) is a certain event, we have that **$p_{i1} + p_{i2} +….. + p_{in} = 1$**, for i=1, 2, …, n.

The row-matrix **$P_k = [p_1^{(k)}\ p_2^{(k)}… p_n^{(k)}]$**, known as the *probability vector* of the Chain, gives the probabilities $p_i^{(k)}$ for the chain to be in state i at step k , for i = 1, 2,…., n and k = 0, 1, 2,…. We obviously have again that **$p_1^{(k)} + p_2^{(k)} + …. + p_n^{(k)} = 1$**.

The following well known Proposition enables one to make *short run forecasts* for the evolution of various situations that can be represented by a finite Markov Chain. The proof of this Proposition is also sketched, just to show to the non expert on the subject reader the strict connection between Markov Chains and Probability (actually Markov chains is considered to be a topic of Probability theory).

*2.1 Proposition:* Under the above notation we have that **$P_{k+1} = P_k A$**, for all non negative integers k.

*Proof:* We shall show that $P_1 = P_0 A$ (1). For this, consider the random variable x = $s_1$, $s_2,……, s_n$, where $s_i$ , i =1,2,…., n are the states of the chain and denote by $F_0, F_1, …$ ,

$F_k$,..... the consecutive steps of the chain. Consider also the events: $E$ = x takes the value $s_1$ in $F_1$, $E_i$ = x takes the value $s_i$ in $F_0$, i = 1, 2, ... , n.

Then $E_i \cap E_j = \emptyset$, i ≠ j, while one (and only one) of the events $E_1, E_2, ..., E_n$ always occurs.. Therefore, by the total probability formula ([7], Chapter 3) we have that $P(E) = \sum_{i=1}^{n} P(Ei) P(E/E_i)$, where $E/E_i$ denote the corresponding conditional probabilities. But $P(E) = p_1^{(1)}$, $P(E_i) = p_i^{(0)}$ and $P(E/E_i) = p_{i1}$. Therefore,
$p_1^{(1)} = p_1^{(0)} p_{11} + p_2^{(0)} p_{21} + .... + p_n^{(0)} p_{n1}$ and in the same way
$p_i^{(1)} = p_1^{(0)} p_{1i} + p_2^{(0)} p_{2i} + ... + p_n^{(0)} p_{ni}$, i = 1, 2,..., n   (2).

Writing the system (2) in matrix form we obtain (1) and working similarly we can show that in general we have $P_{k+1} = P_k A$, for all non negative integers k.

*2.2 Corollary:* Under the same notation we have that **$P_{k+1} = P_0 A^k$** for all non negative integers k, where **$A^0$** denotes the unitary matrix.

*Proof:* Use Proposition 2.1 and apply induction on k.

The following simple application illustrates the above results:

*2.3 Problem:* A company circulates for first time in the market a new product, say K. The market's research shows that the consumers buy on average one such product per week, either K, or a competitive one. It is also expected that 70% of those who buy K they will prefer it again next week, while 20% of those who buy another competitive product they will turn to K next week. Find the market's share for K two weeks after its first circulation, provided that the market's conditions remain unchanged.

*Solution:* We form a Markov Chain having the following two states: $s_1$ = the consumer buys K, and $s_2$ = the consumer buys another competitive product. Then, the transition matrix of the Chain is

$$A = [p_{ij}] = \begin{matrix} & s_1 & s_2 \\ s_1 & \begin{bmatrix} 0.7 & 0.3 \\ 0.2 & 0.8 \end{bmatrix} \\ s_2 & \end{matrix}$$

Further, since K circulates for first time in the market, we have that $P_0 = [0\ 1]$, therefore $P_2 = P_0 A^2 = [0,3\ 0,7]$.
Thus the market's share for K two weeks after its first circulation will be 30%.

## 3. Ergodic Chains

A Markov Chain is said to be an *Ergodic Chain,* if it is possible to go between any two states, not necessarily in one step. It is well known ([4], Chapter 5) that, as the number of its steps tends to infinity (*long run*), an Ergodic Chain tends to an *equilibrium situation*, in which the probability vector $P_k$ takes a constant price **P=[$p_1$ $p_2$ .... $p_n$]**, called the *limiting probability vector* of the chain. Thus, as a direct consequence of Proposition 2.1, the equilibrium situation is characterized by the equality **P = PA**   (3).  Obviously, we have again that **$p_1 + p_2 + ....+ p_n = 1$**.

With the help of Ergodic Chains, one obtains long run forecasts for the evolution of the corresponding phenomena. The following two applications illustrate this situation:

*3.1 Problem:* Reconsidering Problem 2.3, find the market's share for the product K in the long run, i.e. when the consumers' preferences are stabilized.

*Solution:* Obviously the Chain constructed in Problem 2.3 is an Ergodic one. Applying (3) for this chain we get that $[p_1 \ p_2] = \begin{bmatrix} 0.7 & 0.3 \\ 0.2 & 0.8 \end{bmatrix} [p_1 \ p_2]$, which gives $p_1 = 0.7p_1 + 0.2p_2$ and $p_2 = 0.3p_1 + 0.8p_2$, or equivalently that $0.3p_1 - 0.2p_2 = 0$. Solving the linear system of the above equation and of $p_1 + p_2 = 1$ one finds that $p_1 = 0.4$, i.e. the market's share for K in the long run will be 40%.

The next problem concerns the application of a 3-state Ergodic Markov Chain to the production process of an industry:

***3.2 Problem:*** In an industry the production of a certain product is regulated according to its existing stock at the end of the day. Namely, if there exist unsatisfied orders or the stock is zero, then the production of the next day covers the unsatisfied orders plus two more metric units (m.u.). On the contrary, if there exists a non zero stock, there is no production for the next day. We further know that the consumers' demand for the product is either 1 m.u. per day with probability 60%, or 2 m. u. per day with probability 40%. Calculate the probability to have unsatisfied orders in the long run.

*Solution:* Since the maximum product's demand is 2 m. u., the production of the factory at the first day of its function must be 2 m. u. and therefore at the end of the day the stock is either zero or 1 m. u. In the former case the process is repeated in the same way. In the latter case the production of the next day is zero and therefore at the end of this day the stock is either zero (in this case the process is repeated in the same way), or there are unsatisfied orders of 1 m. u. In the last case the production of the next day is 3 m. u., i.e. 1 m. u. to cover the unsatisfied orders of the previous day plus 2 m. u., and so on. It becomes therefore evident that, according to the above rhythm of production, there are three possible situations at the end of each day: $s_1 =$ unsatisfied orders of 1 m.u., $s_2 =$ zero stock and $s_3 =$ stock of 1 m.u. Evidently our problem can be described with an Ergodic Markov Chain having as states the above possible situations $s_i$, $i = 1, 2, 3$. Using the given data it is easy to observe that the transition matrix of the Chain is

$$A = \begin{array}{c} \\ s_1 \\ s_2 \\ s_3 \end{array} \begin{array}{ccc} s_1 & s_2 & s_3 \end{array} \\ \begin{bmatrix} 0 & 0.4 & 0.6 \\ 0 & 0.4 & 0.6 \\ 0.4 & 0.6 & 0 \end{bmatrix}.$$

Let $P = [p_1 \ p_2 \ p_3]$ be the limiting probability vector, then the equality $P = PA$ gives that $p_1 = 0.4p_3$, $p_2 = 0.4p_1 + 0.4p_2 + 0.6p_3$, and $p_3 = 0.6p_1 + 0.6p_2$. Adding the first two of the above equations we find the third one. Solving the linear system of the first two equations and of $p_1 + p_2 + p_3 = 1$ one finds that $p_1 = 0.15$. Therefore the probability to have unsatisfied orders in the long run is 15%.

## 4. Absorbing Markov chains

A state of a Chain is called absorbing if, once entered, it cannot be left. Further a Markov Chain is said to be an *Absorbing Chain* if it has at least one absorbing state and if from every state it is possible to reach an absorbing state, not necessarily in one step.

In case of an Absorbing Chain with k absorbing states, $1 \leq k < n$, we bring its transition matrix A to its *canonical form* $A^*$ by listing the absorbing states first and then we make a partition of $A^*$ of the form $A^* = \begin{bmatrix} I & | & O \\ - & | & - \\ R & | & Q \end{bmatrix}$, where I is the unitary k X k matrix, $O$ is a zero matrix, $R$ is the (n – k) X k transition matrix from the non absorbing to the absorbing states and $Q$ is the (n – k) X (n – k) transition matrix between the non absorbing states.

Denote by $I_{n-k}$ the unitary (n – k) X (n – k) matrix; it can be shown ([4], Chapter 3) that the square matrix $I_{n-k} - Q$ has always a non zero determinant. Then, the *fundamental matrix* of the Absorbing Chain is defined to be the matrix

$$N = [n_{ij}] = (I_{n-k} - Q)^{-1} = \frac{1}{D(I_{n-k} - Q)} adj(I_{n-k} - Q) \quad (4),$$

where $(I_{n-k} - Q)^{-1}$, $D(I_{n-k} - Q)$ and $adj(I_{n-k} - Q)$ denote the inverse, the determinant and the adjoint of the matrix $I_{n-k} - Q$ respectively. We recall that the *adj* $(I_{n-k} - Q)$ is the matrix of the algebraic complements of the transpose matrix of the matrix $I_{n-k} - Q$ ([5], Section 2.4).

It is well known ([4], Chapter 3) that the element $n_{ij}$ of the fundamental matrix N gives the mean number of times in state $s_i$ before the absorption, when starting in state $s_j$ (where $s_i$ and $s_j$ are non absorbing states of the chain).

The above results are illustrated by the next application:

***4.1 Problem:*** An agricultural cooperative applies the following stages for the collection and shelling of a product: $s_1$ = collection, $s_2$ = sorting - refining, $s_3$ = packing and $s_4$ = shelling. The past experience shows that there is a 20% probability that the quality of the collected product is not satisfactory. In this case the collected quantity is abandoned and a new collection is attempted. It is also known that the duration of each of the stages $s_i$, i =1, 2, 3, 4, is on average 10, 4, 3 and 45 days respectively. Find the mean time needed for the completion of the whole process.

*Solution:* The above process can be represented by a finite Markov Chain having as states the stages $s_i$, i=1, 2, 3, 4. This Chain is obviously an Absorbing one and $s_4$ is its unique absorbing state. It is straightforward to check that the transition matrix of the Chain is:

$$A = \begin{array}{c} \\ s_1 \\ s_2 \\ s_3 \\ s_4 \end{array} \begin{array}{cccc} s_1 & s_2 & s_3 & s_4 \\ \begin{bmatrix} 0 & 1 & 0 & 0 \\ 0,2 & 0 & 0,8 & 0 \\ 0 & 0 & 0 & 1 \\ 0 & 0 & 0 & 1 \end{bmatrix} \end{array}, \text{ while its canonical form can be written as}$$

$$A^* = \begin{array}{c} \\ s_4 \\ \\ s_1 \\ s_2 \\ s_3 \end{array} \begin{bmatrix} 1 & | & 0 & 0 & 0 \\ - & | & - & - & - \\ 0 & | & 0 & 1 & 0 \\ 0 & | & 0,2 & 0 & 0,8 \\ 1 & | & 0 & 0 & 0 \end{bmatrix}. \text{ Therefore } I_3 - Q = \begin{bmatrix} 1 & -1 & 0 \\ -0,2 & 1 & -0,8 \\ 0 & 0 & 1 \end{bmatrix} \text{ and (4) gives, after a}$$

straightforward calculation, that $N = \begin{array}{c} s_1 \\ s_2 \\ s_3 \end{array} \begin{bmatrix} 1.25 & 1.25 & 1 \\ 0.25 & 1.25 & 1 \\ 0 & 0 & 1 \end{bmatrix}$ (columns $s_1, s_2, s_3$). Thus, since in this case the Chain is always starting from state $s_1$, the mean number of times in states $s_1$ and $s_2$ before the absorption are 1.25 and in state $s_3$ is 1. Therefore, the mean time needed for the completion of the whole process is $1.25 *(10+4) + 3 + 45 = 65.5$ days.-

When an Absorbing Markov Chain has more than one absorbing states, then the element $b_{ij}$ of the matrix **B = NR = [$b_{ij}$]** gives the probability for the Chain starting in state $s_i$ to be absorbed in state $s_j$ ([4], Chapter 3). This is illustrated in the following example, which is a special case of a "random – walk" problem:

**4.2 Problem:** A supermarket has three storehouses, say $A_1$, $A_2$ and $A_3$ between two cities, say $C_1$ and $C_2$, as it is shown in the below diagram:

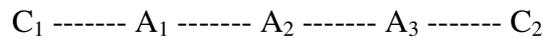

$$C_1 \text{-------} A_1 \text{-------} A_2 \text{-------} A_3 \text{-------} C_2$$

For the delivery of the goods, a truck starts its route every day from one of the storehouses and terminates it to one of the cities. The truck moves each time one place to the right or to the left with the same probability. Find the mean number of stops of the truck to each storehouse during its route and the probability to terminate its route to the city $C_1$, when it starts it from storehouse $A_2$.

*Solution:* We introduce a 5-state Markov Chain having the following states: $s_1$ ($s_5$) = the truck arrives to the city $C_1$ ($C_2$), $s_2$ ($s_3$, $s_4$) = the truck arrives to the storehouse $A_1$ ($A_2$, $A_3$). Obviously the above Chain is an Absorbing one and $s_1$, $s_5$ are its absorbing states. The canonical form of its transition matrix is:

$$A^* = \begin{array}{c} s_1 \\ s_5 \\ \\ s_2 \\ s_3 \\ s_4 \end{array} \begin{bmatrix} 1 & 0 & | & 0 & 0 & 0 \\ 0 & 1 & | & 0 & 0 & 0 \\ - & - & | & - & - & - \\ 0.5 & 0 & | & 0 & 0.5 & 0 \\ 0 & 0 & | & 0{,}5 & 0 & 0.5 \\ 0 & 0.5 & | & 0 & 0.5 & 0 \end{bmatrix}.$$

Then, it is straightforward to check that the fundamental matrix of the Chain is

$$N = (I_3-Q)^{-1} = \begin{array}{c} s_2 \\ s_3 \\ s_4 \end{array} \begin{bmatrix} 1.5 & 1 & 0.5 \\ 1 & 2 & 1 \\ 0.5 & 1 & 1.5 \end{bmatrix}.$$

Thus, since the truck starts its route from the storehouse $A_2$ (state $s_3$), the mean number of its stops to the storehouse $A_1$ (state $s_2$) is 1, to the storehouse $A_2$ (state $s_3$) is 2 and to the storehouse $A_3$ (state $s_4$) is 1.

Further, $B = NR = \begin{bmatrix} 1.5 & 1 & 0.5 \\ 1 & 2 & 1 \\ 0.5 & 1 & 1.5 \end{bmatrix} \begin{bmatrix} 0.5 & 0 \\ 0 & 0 \\ 0 & 0.5 \end{bmatrix} = \begin{matrix} s_2 \\ s_3 \\ s_4 \end{matrix} \begin{bmatrix} \overset{s_1}{0.75} & \overset{s_2}{0.75} \\ 0.5 & 0.5 \\ 0.25 & 0.25 \end{bmatrix}$

Thus the probability for the truck to terminate its route to the city $C_1$ (state $s_1$), when it starts it from store $A_2$ (state $s_3$) is 50%. -

Our last application illustrates the fact that a great care is needed sometimes in order to "translate" correctly the mathematical results of the Markov Chain model in terms of the corresponding real situation.

**4.3 Problem:** In a college the minimal duration of studies is four years. The statistical analysis has shown that there is a 20% probability for each student to be withdrawn due to unsatisfactory performance and a 30% probability to repeat the same year of studies. Find the probability for a student to graduate, the mean time needed for the graduation and the mean time of his/her attendance in each year of studies.

*Solution:* We introduce a finite Markov Chain with the following states: $s_i$ = attendance of the i – th year of studies, i=1, 2, 3, 4 , $s_5$ = withdrawal from the college and $s_6$ = graduation. Obviously the above Chain is an Absorbing Markov Chain, and $s_5$ and $s_6$ are its absorbing states. The canonical form of its transition matrix is

$$A^* = \begin{matrix} s_5 \\ s_6 \\ \\ s_1 \\ s_2 \\ s_3 \\ s_4 \end{matrix} \begin{bmatrix} \overset{s_5}{1} & \overset{s_6}{0} & | & \overset{s_1}{0} & \overset{s_2}{0} & \overset{s_3}{0} & \overset{s_4}{0} \\ 0 & 1 & | & 0 & 0 & 0 & 0 \\ - & - & | & - & - & - & - \\ 0.2 & 0 & | & 0.3 & 0.5 & 0 & 0 \\ 0.2 & 0 & | & 0 & 0.3 & 0.5 & 0 \\ 0.2 & 0 & | & 0 & 0 & 0.3 & 0.5 \\ 0.2 & 0 & | & 0 & 0 & 0 & 0.3 \end{bmatrix}$$

Further, using a PC mathematical package to make the necessary calculations quicker, it is straightforward to check that the fundamental matrix of the chain is

$N = (I_4 - Q)^{-1} = \begin{matrix} s_1 \\ s_2 \\ s_3 \\ s_4 \end{matrix} \begin{bmatrix} \overset{s_1}{1.429} & \overset{s_2}{1.02} & \overset{s_3}{0.729} & \overset{s_4}{0.521} \\ 0 & 1.429 & 1.02 & 0.729 \\ 0 & 0 & 1.429 & 1.02 \\ 0 & 0 & 0 & 1.429 \end{bmatrix}$ and therefore to find that

$B = NR = \begin{matrix} s_1 \\ s_2 \\ s_3 \\ s_4 \end{matrix} \begin{bmatrix} \overset{s_5}{0,74} & \overset{s_6}{0,261} \\ 0,636 & 0,365 \\ 0,49 & 0,51 \\ 0,286 & 0,715 \end{bmatrix}$.

Observing the fundamental matrix N of the chain, one finds that $n_{13}$=0.729 and $n_{14}$=0.521, i.e. for a first year student of the college the mean time of attendance in the third and fourth year of studies is less than one year! However this is not embarrassing, because there is always a possibility for a student to be withdrawn from the college due to unsatisfactory performance before entering the third, or fourth, year of studies.

Since $n_{11} = n_{22} = n_{33} = n_{44} = 1.429$, we find that the mean time of attendance of a student in each year of studies is 1.429 years, while the mean time needed for his/her graduation is 1.429 * 4 = 5.716 years.

Further, observing the matrix B one finds that $b_{15} = 0.74$, i.e. the probability of a student to graduate is 74%.

## 5. A Brief Account of our Relevant Researches

There are very many applications of Markov Chains reported in the literature for the solution of real world problems in almost every sector of the human activities. However, a complete reference to all, or at least to the most important of them, is out of the scope of the present paper. Here we shall only restrict to a brief presentation of the most important of our personal research results concerning applications of Finite Markov Chains to *Management*, *Mathematical Education* and *Artificial Intelligence*.

In fact, in Voskoglou and Perdikaris [10, 11], and in Perdikaris and Voskoglou [6] the *Problem Solving* process is described by introducing an Absorbing Markov Chain on its main steps. Also, in Voskoglou [13] an Ergodic Chain is introduced for the study of the *Analogical Reasoning* process in the classroom through the solution of suitably chosen problems, while in Voskoglou [14] an Absorbing Markov Chain is utilized for the description of the process of *Learning a Subject Matter* in the classroom.

Further, in Voskoglou [12] an Absorbing Markov Chain is introduced to the major steps of the *Modelling* process applied for the study of a real system. An alternative form of the above model is introduced in Voskoglou [16] for the description of the *Mathematical Modelling* process in the classroom, when the teacher gives such kind of problems for solution to students. In this case it is assumed that the completion of the solution of each problem is followed by a new problem given to the class, which means that the modelling process restarts again from the beginning. Under this assumption the resulting Markov Chain is an Ergodic one.

In Voskoglou [15, 18] an Absorbing Markov chain is introduced to the main steps of the *Decision Making* process and examples are presented to illustrate the applicability of the constructed model to decision making situations of our day to day life. .Finally in [17] the theory of Markov Chains is used for the description of the *Case -Based Reasoning (CBR)* process, which is applied for the solution of a new problem in terms of already known from the past solutions of similar problems. A measure is also obtained in [17] for testing the efficiency of a CBR system designed to work with the help of computers.

## 5. Final conclusions

The following conclusions can be obtained from the discussion performed in this review paper concerning applications of finite Markov Chains to Management:

- The theory of Markov Chains is a successful combination of Linear Algebra and Probability theory, which enables one to make short and long run forecasts for the evolution of various phenomena of the real world characterized by only "one step memory".
- The *short run* forecasts are obtained, regardless the type of the Chain, by calculating its transition matrix and the probability vector of its corresponding step (basic form of the Markov Chain model). On the contrary, the *long run* forecasts (equilibrium situation of the Chain) are obtained, in case of Ergodic Chains only, by calculating their limiting probability vectors.
- In case of an *Absorbing Markov Chain* one proceeds to the study of the corresponding problem by forming the canonical form of its transition matrix and calculating the fundamental matrix of the chain, the entries of which give the mean number of times in each non absorbing state before absorption, for each possible non absorbing starting state. When an absorbing chain has more than one absorbing states, then the transition matrix from the non absorbing to the absorbing states enables us to calculate the probabilities for the chain to reach a certain absorbing state, when it starts from a certain non absorbing state.